# A Conversation with Ulf Grenander

**Nitis Mukhopadhyay**

*Abstract.* Ulf Grenander was born in Vastervik, Sweden, on July 23, 1923. He started his undergraduate education at Uppsala University, and earned his B.A. degree in 1946, the Fil. Lic. degree in 1948 and the Fil. Dr. degree in 1950, all from the University of Stockholm. His Ph.D. thesis advisor was Harald Cramér. Professor Grenander is well known for pathbreaking research in a number of areas including pattern theory, computer vision, inference in stochastic processes, probabilities on algebraic structures and actuarial mathematics. He has published more than one dozen influential books, of which *Statistical Analysis of Stationary Time Series* (1957, coauthored with M. Rosenblatt), *Probabilities on Algebraic Structures* (1963; also in Russian) and *Abstract Inference* (1981b) are regarded as classics. His three-volume lecture notes, namely, *Pattern Synthesis* (vol. I, 1976), *Pattern Analysis* (vol. II, 1978) and *Regular Structures* (vol. III, 1981a; also in Russian) created and nurtured a brand new area of research. During 1951–1966, Professor Grenander's career path took him to the University of Chicago (1951–1952), the University of California–Berkeley (1952–1953), the University of Stockholm (1953–1957), Brown University (1957–1958) and the Institute for Insurance Mathematics and Mathematical Statistics (1958–1966) as its Professor and Director. From 1966 until his retirement he was L. Herbert Ballou University Professor at Brown University. Professor Grenander also held the position of Scientific Director (1971–1973) of the Swedish Institute of Applied Mathematics. He has earned many honors and awards, including Arhennius Fellow (1948), Fellow of the Institute of Mathematical Statistics (1953), Prize of the Nordic Actuaries (1961), Arnberger Prize of the Royal Swedish Academy of Science (1962), Member of the Royal Swedish Academy of Science (1965), Guggenheim Fellowship (1979) and Honorary Fellow of the Royal Statistical Society, London (1989). He has delivered numerous prestigious lectures, including the Rietz Lecture (1985), the Wald Lectures (1995) and the Mahalanobis Lecture (2004). Professor Grenander received an Honorary D.Sc. degree (1993) from the University of Chicago and is a Fellow of both the American Academy of Arts and Sciences (1995) and the National Academy of Sciences, U.S.A. (1998). Professor Grenander's career, life, passion and hobbies can all be summarized by one simple word: Mathematics.

*Nitis Mukhopadhyay is Professor, Department of Statistics, University of Connecticut, Storrs, Connecticut 06269-4120, USA e-mail: mukhop@uconnvm.uconn.edu. Ulf Grenander's e-mail address at Brown University is ulf-grenander@cox.net.*

The following conversation took place on May 11 and 12, 2002 at Professor Ulf and Mrs. Paj Grenander's lovely home in Providence, Rhode Island, only minutes from the Brown University campus.

## UPBRINGING: VASTERVIK, SWEDEN

**Mukhopadhyay:** Ulf, shall we start at the very beginning? When were you born and what was your birthplace like?

**Grenander:** I was born on July 23, 1923 in a small town called Vastervik, situated on the east coast of Sweden. It was a very lovely place actually, but there was nothing important or striking about the town itself. We still keep our summer house a little outside of that place and we continue to enjoy it very much.

**Mukhopadhyay:** Would you mention your parents and some backgrounds?

**Grenander:** My father, Sven, was a very interesting person. His background was in mathematics and physics, and he was trained as a meteorologist. Later in life, he had only limited interest in these subjects even though he had a Fil. Dr. degree in meteorology from the University of Stockholm. A Fil. Dr. degree recipient from University of Stockholm is awarded a doctoral ring along with the diploma. I always wear my father's doctoral ring. It has the inscription 1911 inside, the year of my father's graduation, and also has the inscription 1950, my graduating year.

**Mukhopadhyay:** Did your father change his field?

**Grenander:** In a way, he did. His great passion of life was actually sailing. Sailing was and still is an expensive hobby. So, my father had to find a clever way to finance this hobby.

**Mukhopadhyay:** How did he manage to do that?

**Grenander:** This was a big boat, a German built yawl Senta, and he had a crew of ten young boys. Instead of paying the boys, my father made those teenagers pay him for the privilege of being on the boat in the first place. Those teenagers were ready to accept the adventure as a challenge and their parents had to pay my father up front.

**Mukhopadhyay:** Do you recall interactions with your father?

**Grenander:** I did not have much mathematical interactions with my father. My father was a school teacher of physics. We had many interesting conversations on subjects other than mathematics. I learned a lot from him about history, politics and social science. One day he enthusiastically told me about Mendelian theory. This learned man was keenly interested in things like this, but he was not interested in mathematics for its own sake. As far as mathematics goes, I was on my own.

**Mukhopadhyay:** How about your mother?

**Grenander:** My mother, Maria, was a housewife. She did not go outside to work toward a career. At that time, women did not go out much to earn a living.

**Mukhopadhyay:** Did you have brothers or sisters?

**Grenander:** I had one brother, Nils. Since he was six years older than I was, I did not have many opportunities to interact with him during my formative years. He became a lawyer specializing in laws regulating the shipping industry. This used to be an esoteric field for lawyers to pursue.

**Mukhopadhyay:** How about your early schooling?

**Grenander:** I attended a local school where I went through all the grades from K through 12. Practically every Scandinavian school was state-operated and my school was no exception.

**Mukhopadhyay:** What did you study in higher grades?

**Grenander:** In higher grades, I pursued the classical curriculum, which meant that it included mainly languages. I mean I had to keep learning one language after the other! Nitis, you being from Indian

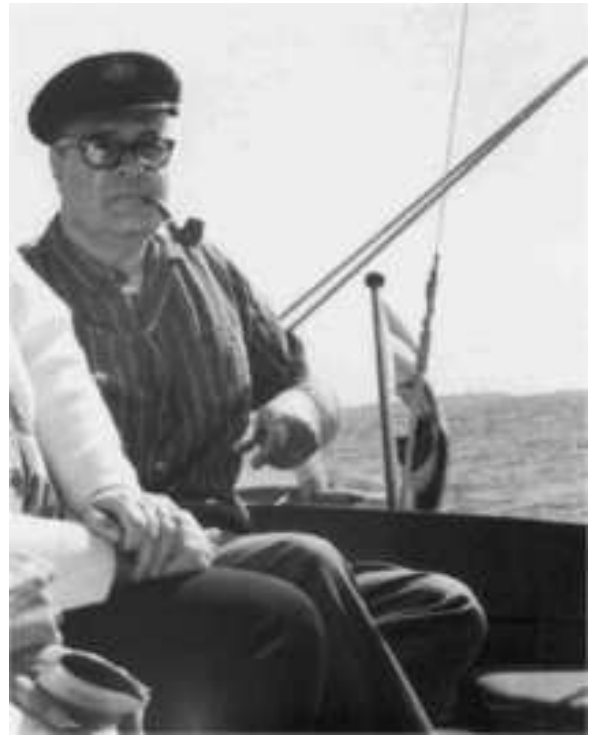

Fig. 1.   *Ulf Grenander sailing in the Baltic.*

origin, you must know several languages. Do you know Sanskrit?

**Mukhopadhyay:** Yes, I had to learn Sanskrit. I studied Bengali, my mother tongue, and also Hindi. I learned English from very early childhood. Ulf, what languages did you learn in school?

**Grenander:** Of course, there was much emphasis on Scandinavian languages. There were three of them. Then, I started learning German, my first foreign language, and I studied it for eight years. Next came English and French. Latin was very important and I rather liked it because of the logical nature of its grammar. I returned to the study of grammar much later though. Unfortunately, there was hardly any mathematics in the curriculum.

**Mukhopadhyay:** Did you not learn any mathematics in high school?

**Grenander:** I mean that I did not receive any rigorous early schooling in mathematics. We studied Euclid's geometry, but I did not understand the material and I did not see the point.

**Mukhopadhyay:** Did its axiomatic approach bother you?

**Grenander:** I did not really appreciate why Euclid had built his theory the way he did. I asked myself, "Why must one build a machinery like that?" My appreciation of Euclidean geometry came much later.

**Mukhopadhyay:** Did you study science in school?

**Grenander:** I did not learn much of science either. There was a little bit of physics in the curriculum, but we had nothing in chemistry or biology.

## EXPERIENCING MATHEMATICS AND STATISTICAL MECHANICS

**Mukhopadhyay:** It appears that neither school nor immediate family helped much in nurturing your mathematical upbringing. What triggered the urge to explore the wide world of mathematics?

**Grenander:** Like most boys, I used to experiment with electricity. Once I had built a small radio transmitter that actually transmitted a message, something that was illegal to do! When this episode became known, my father was not too happy about it, but he was also feeling a little guilty because he was the one who had helped me to gather some of the essential components that went inside my transmitter. (Laughs)

I picked up mathematics on my own. My school library's collections mainly catered to the classical disciplines, but strangely it also housed an enormous mathematical encyclopedia (*Encyklopedie der Mathematischen Wissenschaften*), published in both German and French. These covered materials from the period 1900 through 1939, I believe. I saw some magnificent things there!

**Mukhopadhyay:** What were some of those magnificent things?

**Grenander:** When I was in the tenth grade, I read a great article by Ehrenfest and Ehrenfest, a married couple. This article gave the first rigorous treatment of statistical mechanics and it was amazing. It gave a very abstract treatment of gas particles consisting of small squares moving in a plane, simply going up or down and nothing else, with a uniform velocity but colliding! Starting with a very abstract formulation, Ehrenfest and Ehrenfest derived the laws of thermodynamics. I thought that it was wonderful. It was remarkable how statistical mechanics could be founded on just few general principles.

**Mukhopadhyay:** Did you pursue works of Fermi, Dirac or Boltzmann?

**Grenander:** A year later, I certainly familiarized myself with some of Boltzmann's ideas. After grade eleven, that is in 1940–1941, I started to pick up quantum mechanics. I managed to buy some enormous German textbooks. Because of the war, none was available in English or French. Of course, there was no such material in Swedish. Statistical mechanics was definitely my first love.

I recall that this was the first time I came across the ideas of probability and statistics. Of course in those days, at least within the mathematical circles, the word "statistics" meant probability theory and not inference.

Graduating from high school in 1942 was serious business at the time. This was similar to Baccalauria from France and it was not easy.

## ARNE BEURLING: UPPSALA UNIVERSITY

**Mukhopadhyay:** Next you enrolled in the undergraduate program at Uppsala University. What do you recall about the transition from high school to college?

**Grenander:** I joined Uppsala University in 1942 during World War II. The war was too close to home in Northern Europe and the situation was very serious. Everyone was nervous and felt unsure about the future.

Initially, I decided to study mathematics and mechanics, but around 1944–1945, I changed my emphasis to mathematics, statistics and probability theory. At that time, Stockholm University was the only Scandinavian academic institution offering undergraduate education in mathematical statistics.

**Mukhopadhyay:** Surely, there were many first-rate mathematicians there at the time. Who had influenced you the most?

**Grenander:** I had a wonderful experience. I came in contact with the mathematics professor, Arne Beurling. He was a great analyst. He became one of the best code breakers during World War II by breaking the German strategic code (warning for Barbarossa to British intelligence). At that time, I did not realize what a great scientist Beurling was. He offered a graduate seminar and what a joy it was when he allowed me to attend the seminar even though I was an undergraduate student. There were only four of us attending Beurling's lectures. Later on, one of the students in attendance became very well known in probability theory. This was G. Esseen of the famous Berry–Esseen theorem!

**Mukhopadhyay:** What did you find so special about Beurling's lectures?

**Grenander:** This small seminar group met at eight o'clock in the evening, but only once every two weeks, spanning a couple of hours at a time. Most of the time Beurling spoke himself about his own recent research. The lectures were absolutely exciting. One could see a brilliant mind at work. I also attended Beurling's other lectures on topics that later became the field of spectral synthesis.

I was indeed fortunate and lucky to have met Beurling that early. His radiating intellectual power influenced me like no one else's.

## HARALD CRAMÉR'S INSTITUTE AND MILITARY SERVICE

**Mukhopadhyay:** Did you finish your undergraduate degree at Uppsala?

**Grenander:** I finished the first two years in Uppsala and then transferred over to Harald Cramér's Institute for Insurance Mathematics and Mathematical Statistics in Stockholm in 1944–1945. This institute was famous for its programs. After I spent two years in this institute, I received my undergraduate degree from the University of Stockholm.

**Mukhopadhyay:** Other than Harald Cramér, who else comes to your mind?

**Grenander:** I recall an experience I had with Harald Bohr, a brother of the famous Nobel Laureate Niels Bohr. I was told that Harald Bohr was actually the brighter of the two and this mathematician discovered the theory of periodic functions. Harald Bohr was of Jewish origin and lived in Denmark, which was occupied by Germany, so he had to flee. As a matter of fact, all the Danish Jews were rescued in one big operation overnight and brought to Sweden.

Harald Bohr came to visit Cramér's institute when I was there as an undergraduate student. He gave lectures on his own work, and that was the first time I heard anyone talking about functions as points in a space. I remember especially one of his lectures on almost periodic functions. You know that the French painter Matisse once said, "the purpose of an artist is to decorate the surface." At the time, Bohr did not have PowerPoint at his disposal, but he had colored chalk. He would start writing at 10:15 in the upper left corner of the chalkboard and then using all sorts of possible colors he would draw his functions and other mathematical structures to fill the whole board. Exactly at 11:00, he would reach the lower right corner of the chalkboard and the lecture would come to an end. The way he talked about what is now called functional analysis was fascinating.

**Mukhopadhyay:** Was harmonic analysis tied in too?

**Grenander:** Yes, it certainly was. I saw plenty of Fourier analysis, but it was mainly through Norbert Wiener's work on generalized harmonic analysis.

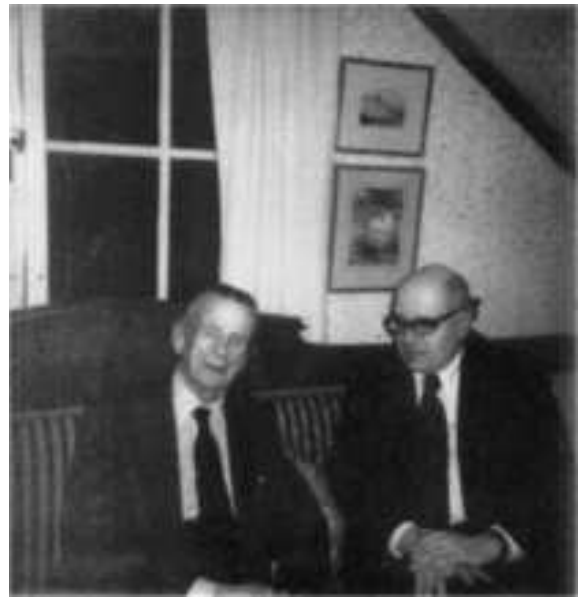

FIG. 2. *Harald Cramér and Ulf Grenander, 1984.*

**Mukhopadhyay:** Ulf, please excuse me for asking this. You were a young undergraduate at the time. You certainly appreciated the daunting beauty in Bohr's lectures, but did you feel that you were mathematically equipped to grasp the depth of those lectures?

**Grenander:** By then, I was mathematically well prepared. The most difficult part of those lectures was that Bohr talked in Danish, which I could understand perhaps, but then he was using the Copenhagen dialect. The Copenhagen accent and dialect were much more difficult to understand than the actual mathematics. (Laughs)

**Mukhopadhyay:** What was the relationship between Cramér's Institute and the University of Stockholm?

**Grenander:** The Institute was actually a department or a center within the university and it was initiated by the donations received from insurance companies.

Cramér was certainly the only senior faculty member in mathematical statistics in Stockholm and for that matter he was the only senior one in this field in all of Scandinavia. This was 1944–1945 and his big book, *Mathematical Methods of Statistics* (Cramér, 1946), came out little later. Perhaps he was still working on his book. Of course, I went to Cramér's lectures and they were wonderful.

**Mukhopadhyay:** What did he teach?

**Grenander:** He taught very much the material that was in his book. This is all standard material now. He introduced plenty of large-sample theory, estimation theory and tests of linear hypotheses including the analysis of variance. I recall, however, that there was no mention of time series analysis or nonparametric statistics.

**Mukhopadhyay:** Cramér's book included everything that mathematical statistics had to offer at the time. It was state-of-the-art presentation, was it not?

**Grenander:** Nitis, actually Cramér's lectures were always like that.

I should mention that I served in the army right after I finished my undergraduate degree from Stockholm.

**Mukhopadhyay:** How did you cope with the military assignment?

**Grenander:** Military service was not so bad! I discovered the pleasures of outdoor life instead of being a nerd. It was fun to shoot cannons, especially since no one was shooting back at me! After serving about a year in 1945, I came back to Stockholm and got myself enrolled for graduate studies. But, unfortunately for me, around that time, Cramér became the President of the University of Stockholm and then he hardly had any time to spare for advising students.

## ACTUARY, INSURANCE AND INFLUENCE OF CRAMÉR

**Mukhopadhyay:** Harald Cramér started out as an actuary and he had a big influence on you, did he not?

**Grenander:** Oh yes, he certainly had a big influence on me. Yes, Cramér was an actuary. Perhaps I should tell you an interesting story. Young Cramér was working on analytic number theory, and he and Harald Bohr had written a very influential article in the mathematical encyclopedia that I mentioned earlier. This was a beautiful paper and Cramér would have continued with analytic number theory to become an assistant professor and so on. This career path was almost set, but a problem erupted between Cramér and his teacher, Mittag-Leffler. Apparently, Cramér had published a paper on number theory as his own work and Mittag-Leffler claimed that it was not quite original. The episode was apparently very embarrassing. Mittag-Leffler supposedly told Cramér, "I am going to see that you will never get any academic position in Northern Europe."

**Mukhopadhyay:** (Laughs) What did young Cramér do then?

**Grenander:** Well, in analytic number theory, one counts and derives frequency ratios, and the whole approach is very much like probability theory. They are not the same, but some of the tools, for example, Fourier analysis, are similar. He found a natural opening in a field where probability theory was used, namely, insurance, and he became an insurance actuary in a company. He never left that. For all those years he held that position, but he did so only as a part timer in later life. Cramér used to say, "That was his great luck," because this was where he discovered probability theory and realized how useful it could be.

**Mukhopadhyay:** Were you drawn into some of the things that he was doing at the time?

**Grenander:** Well, he was interested in modeling demographics and in the construction of life tables, but not in a pedestrian way of determining premiums and such for insurance purposes. He had to do some of that and he taught me the ropes. You see, I was initially trained as an actuary too! I actually practiced as a consulting actuary for a short while! But Cramér introduced the idea of what was known as the "ruin problem."

**Mukhopadhyay:** How would a layman connect a "ruin problem" with insurance strategies?

**Grenander:** An insurance company loses money when it pays off an insured customer, but gains money from premiums and interests on its reserve. But it does not want to lose too much cash to cause near bankruptcy. Some people with life insurance die early, whereas some do not. Some existing customers cancel policies, whereas new customers buy new policies. An executive may want to know, for example, the probability of the company going bankrupt or that its cash reserve will go under a certain threshold in one single year. Of course such probabilities should be small.

**Mukhopadhyay:** This is a very dynamic stochastic system.

**Grenander:** Yes, it is. Actually, Cramér taught us to think of a random process governed by differential equations. It was called Thiele's differential equations which are now referred to as stochastic differential equations. The central limit theorem does not work well because one is now interested in approximating the outer extreme or the tail of a probability distribution. Cramér invented large deviation theory and this was exactly the right tool that was needed in the insurance industry. It was a great beginning. A few years before Cramér died, he used to come and visit when I was at the Mittag-Leffler Institute outside of Stockholm. I once mentioned to Cramér that the field of large deviations had turned into a minor industry, but he did not seem to believe it. He probably thought that I was just being polite.

**Mukhopadhyay:** Was Cramér's institute a one-man show in some sense?

**Grenander:** Yes, this was almost a one-man show. Cramér was surrounded by students and occasionally there were younger colleagues too. The institute had a visitor there from 1936 to 1940, and this was William Feller. Feller came from Germany and Cramér took him under his wings. I think it was Cramér who told him that he should not continue working on differential geometry, which he did before. Cramér advised him to move to probability theory, the field of the future, and then Feller started working on probability theory.

I did not overlap with Feller at that time, but later I had many contacts with him. He used to write to me in Swedish! (Laughs)

## PH.D. DEGREE FROM STOCKHOLM

**Mukhopadhyay:** Did you start working with Cramér on your thesis?

**Grenander:** I did not work with him because you would recall that at that time he was the President of the University of Stockholm and he was completely consumed by administrative work. We had very little adviser–student interaction, but he was always supporting me.

**Mukhopadhyay:** I recall that later Cramér became the Chancellor of the Swedish system of higher academic institutions.

**Grenander:** Exactly, but then in terms of availability of Cramér, an already bad situation became much worse!

He was the President of the University when I started, and I did not see him very much. That was bad enough for me. I wrote a Licentiate thesis, something like a Ph.D. thesis I suppose, in Swedish. I am glad I did not write it in a language that anyone else could read it easily. (Laughs)

**Mukhopadhyay:** Why was that?

**Grenander:** I worked by myself for two years on my Licentiate thesis that developed integration on abstract spaces, such as a Banach space or a lattice. It was so abstract! I am almost embarrassed when I read it myself now. Cramér signed the final copy of the thesis (1948).

You see that the second half of my Licentiate thesis was more interesting where I started to think about ways to make statistical inference for stochastic processes. This preliminary investigation led to my Ph.D. thesis, Stochastic Processes and Statistical Inference (Grenander, 1950). I continued researching in this area for nearly two more years after that.

**Mukhopadhyay:** At some point, you interacted with Karhunen, a Finnish mathematician.

**Grenander:** There was an experienced Finnish individual who was much older than me. This was Karhunen, well-known for the Karhunen–Loève expansion of stochastic processes. He visited Cramér's institute for a year from the University of Helsinki. We met regularly to discuss mathematics. He was extremely good and I really appreciated his ideas. Unfortunately, there was no vacancy for a mathematical statistician, and so Karhunen also went to work for insurance. Actually, he became the CEO of the largest insurance company in Finland. Later, I once asked him, "How can you stand the kind of work that is demanded of you?" Karhunen quickly replied, "Nowadays, I solve problems, but that was exactly what I used to do 20 years ago! The difference is that I solve different kinds of problems now." (Laughs)

## UNIVERSITY OF STOCKHOLM FACULTY: KOLMOGOROV'S APPRECIATION OF PH.D. THESIS

**Mukhopadhyay:** What did you do after finishing your Ph.D. degree?

**Grenander:** After receiving the Ph.D. degree, I became an assistant professor for a year at the University of Stockholm. In Scandinavia, at that time, assistant professors used to carry a light teaching load amounting to two hours of teaching per week.

**Mukhopadhyay:** Ulf, what did you teach?

**Grenander:** I taught my own stuff, namely, inference for stochastic processes, to the graduate students, and there were very few of them around. The custom was very different then. An assistant professor was expected to teach graduate students the advanced material from their own research, whereas a full professor would teach the undergraduates carrying a heavier load, perhaps four hours per week.

I may mention that initially my thesis did not draw any real attention until Kolmogorov got hold of it. Later, the field grew and my thesis was appreciated more. Rao's (2000) recent book, *Stochastic Processes*: *Inference Theory*, is based on my thesis.

**Mukhopadhyay:** How did Kolmogorov come across your Ph.D. thesis?

**Grenander:** Kolmogorov had contacts with Cramér and he visited Stockholm often, so I came to know him well. That was a fantastic experience.

Kolmogorov did not give many lectures. I found the ones that I attended very difficult. I remember especially one lecture that he gave much later on functions being simple or complex. At the time, I could not understand what he meant by that. It material actually led to wonderful work that became known as Kolmogorov's complexity theory. This was a very fundamental piece of work.

**Mukhopadhyay:** Did you have any opportunity to discuss technical matters with Kolmogorov? Did he offer any advice?

**Grenander:** Let me start by saying yes to both questions. Kolmogorov was very encouraging and he was the first person who told me to continue my work on inference for stochastic processes. Later, I discovered that he traveled in Eastern Europe talking to other people about my Ph.D. thesis. That is how my Ph.D. thesis became known.

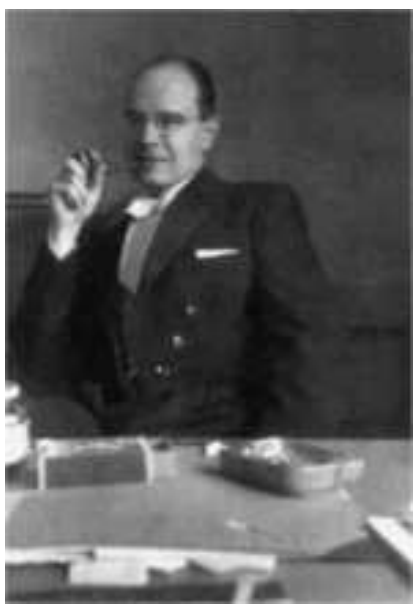

FIG. 3. *Ulf Grenander defending his Ph.D. thesis, Stockholm, 1950.*

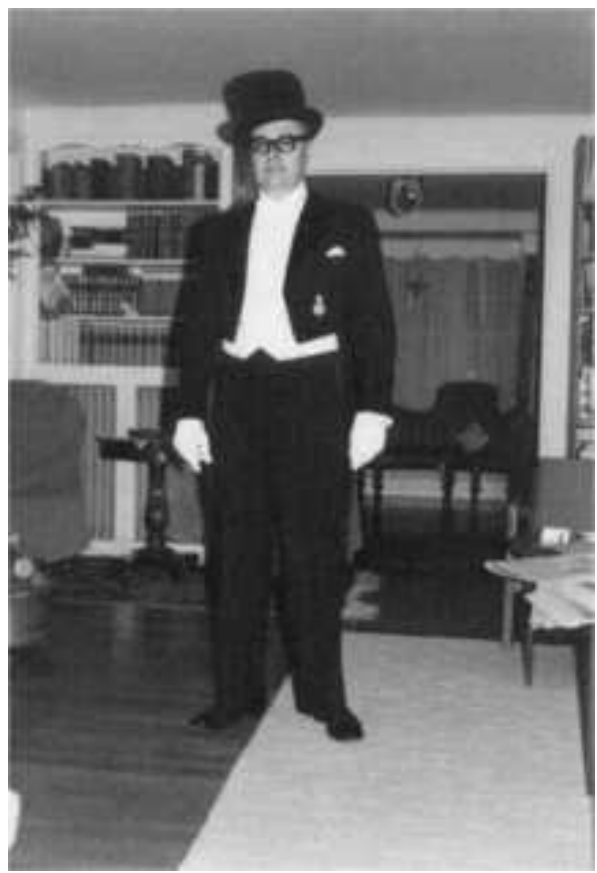

FIG. 4. *Ulf Grenander as a new professor at Stockholm, 1959.*

**Mukhopadhyay:** Was language a barrier in exchanging ideas with Kolmogorov?

**Grenander:** It might have been somewhat of a barrier. When I met Kolmogorov, we talked in German. He could read English but he pronounced in French and that would sound rather strange.

Kolmogorov used to come to Stockholm by boat because he did not like to fly. One time in the 1960s while he was in Stockholm, he told me that he enjoyed the visit very much except for one thing. He asked, "Why don't your young people come and talk to me?" I tried to make him feel comfortable by explaining, "Well, you are a great scientist. These younger people are just students and hence they seem to shy away because of your fame and stature, I suppose." Kolmogorov immediately snapped and said, "That is the most ridiculous thing I have ever heard." A little later I was walking him to the boat to go back to Leningrad. On our way, we were discussing his time, early 1930 or perhaps even earlier, in Gottingen. I said, "This was a wonderful group of mathematicians, the best in the world I suppose. You must have enjoyed the company of David Hilbert." Kolmogorov replied immediately without a second thought, "Now, don't you realize that this was the great one and only David Hilbert, and I was just a student? I could not have interacted with Hilbert!" (Laughs)

## VISITS TO CHICAGO AND BERKELEY (1951–1953): ENCOUNTERS WITH NEYMAN, ROSENBLATT AND SZEGÖ

**Mukhopadhyay:** So, Kolmogorov's personal attention gave you a big break.

**Grenander:** It surely did. Kolmogorov obviously appreciated my Ph.D. thesis. Right around that time, people in America found out about my work as well, and I was asked to come and visit the University of Chicago. That I believe was due to the generosities of both Kolmogorov and Neyman. I came to Chicago as a visiting assistant professor for a year in 1951–1952. I had a wonderful time in the Department of Statistics, but it did not bear that name at the time. It was called the Committee on Statistics, which was headed by Allen Wallis. This group included Bill Kruskal, Murray Rosenblatt, Charles Stein, Leo Goodman—about ten young individuals in all. This was a very stimulating group and I loved it. Charles Stein and I shared an office in Chicago. It gave me an opportunity to share many ideas with him. I was so impressed with his depth of understanding. Joe Hodges from Berkeley was also visiting Chicago and we became good friends.

**Mukhopadhyay:** Did you like renting a place from Jimmy Savage?

**Grenander:** In Chicago, I rented an apartment from Jimmy Savage and so I had the rare opportunity to read mathematical books from his personal library. I found Pontryagin's *Topological Groups*, which opened a new world for me. It really did. It is a good idea to live in another person's house if it also has a good library!

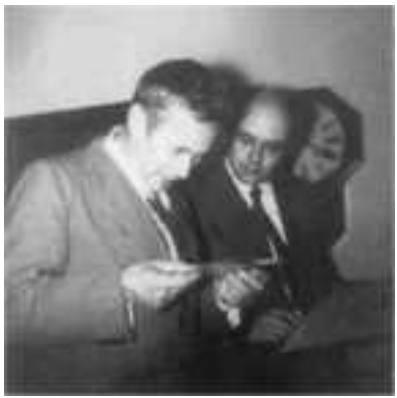

Fig. 5. *A. N. Kolmogorov (left) and Ulf Grenander.*

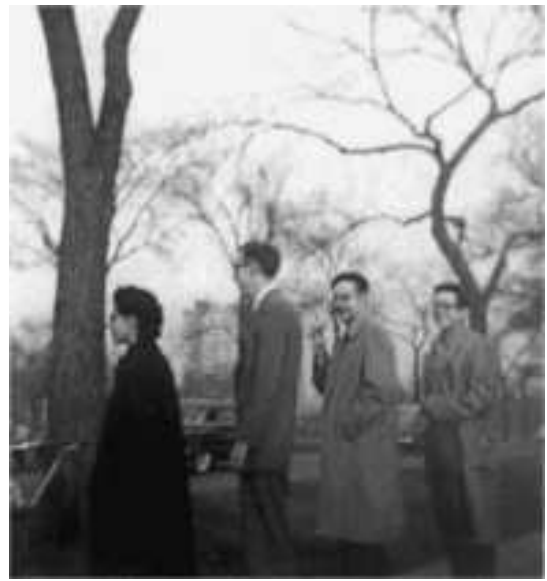

Fig. 6. *Forefront (from left to right): Ady Rosenblatt, Charles Stein, Ulf Grenander and Murray Rosenblatt, at Chicago, 1951.*

**Mukhopadhyay:** You had collaborated with Murray Rosenblatt. What urged you to work in time series?

**Grenander:** My own interest in time series began after reading the wonderful papers of C. E. Shannon and Stuart Rice in the *Bell System Technical Journal*. They dealt with information and how information moved from one channel to another. This and my previous interest in stochastic processes made it quite natural for me to work in time series at that time.

**Mukhopadhyay:** I understand that you had a contact with Jerzy Neyman. Please tell me about it.

**Grenander:** Neyman had seen some work that I had been doing, not involving stochastic processes, but rather my early material on time series. I was beginning my research on time series using the work of G. Szegö, namely the Toeplitz forms. When I was in Chicago, I received a call one day from Neyman to go and spend a year at Berkeley, which of course I did. I spent the year 1952–1953 as a visiting associate professor at the University of California–Berkeley, and that visit was also very stimulating.

**Mukhopadhyay:** In Berkeley, who were some of the people that you met?

**Grenander:** Other than Neyman, I saw Erich Lehmann, Charles Stein, Michel Loève, Joe Hodges and many others. It was a large group of people. At that time, Berkeley and Stanford used to hold joint seminars, perhaps once or twice a month. I came to meet Szegö during such a joint seminar.

Szegö was the chairman of the Mathematics Department at Stanford. He was working on Toeplitz forms since 1923 and I was using them since 1951. He expressed interest in writing a book together on this topic. We wrote the book *Toeplitz Forms and Their Applications* in 1952–1953, but it came out several years later in 1958 (Grenander and Szegö, 1958).

**Mukhopadhyay:** What was Szegö's personality like?

**Grenander:** Well, it was very hard work to finish the book. Szegö was very demanding and extremely precise, but elegant. He had sharp analytical intelligence. He did not favor unnecessary complications. He stuck to his own views. He was a typical representative of Central European scholarship, a cultural phenomenon that does not really exist today.

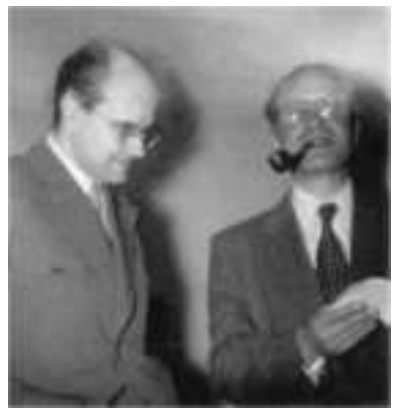

Fig. 7. *Ulf Grenander and M. Loève (right) at Berkeley, 1952.*

**Mukhopadhyay:** What else about that time in Berkeley?

**Grenander:** You will recall that I interacted with Murray Rosenblatt in Chicago and jointly wrote the book *Statistical Analysis of Stationary Time Series*. This book also came out much later, in 1957 (Grenander and Rosenblatt, 1957).

Neyman encouraged me to pursue any area that I personally enjoyed to work in. It was characteristic of Neyman to encourage the younger generation. It is the one thing that first comes to my mind when I think of him.

## RETURN TO UNIVERSITY OF STOCKHOLM (1953): CONSULTING ACTUARY

**Mukhopadhyay:** What happened upon your return to Stockholm in 1953?

**Grenander:** Rosenblatt came to visit Stockholm and stayed for a year. That must have been 1954, I suppose. During 1955–1957, I was an assistant professor at the University of Stockholm and I was also a consulting actuary, which I found extremely interesting.

**Mukhopadhyay:** What was so interesting about consulting as an actuary?

**Grenander:** By that time, Cramér became the Chancellor and I saw him even less frequently. The insurance company I was involved with asked me to look into the construction of mortality tables. Apparently the existing tables were becoming obsolete. At that time everything was based on a three-parameter representation called Makeham's formula that worked to some extent, but there was no theory behind that formula. It was completely empirical in nature.

I started thinking whether I could do away with any analytic expression, that is, I was searching for a suitable nonparametric statistical procedure. We knew that the mortality intensity was an increasing function of age. I asked myself, "What is the best estimator of mortality intensity in the sense of maximum likelihood?" This ultimately led to estimation of unimodal densities (Grenander, 1956, 1957).

**Mukhopadhyay:** Did you work on problems arising in automobile insurance?

**Grenander:** Yes, I became interested in problems related to automobile insurance. Some companies used what was known as a bonus system. If one reported an accident to the company, soon his premium went up. If there was no reported accident in a year, the next premium went down slightly. I believe that the automobile insurance industry in Sweden still uses a system with seven classes or states. It is obviously a Markov chain with seven states and one could ask whether one should report an accident. One has to report an accident to the police, but not necessarily to his insurance carrier. This amounts to a problem on optimization that I solved (Grenander, 1958).

**Mukhopadhyay:** Will it be fair to say that you like to develop crucial mathematical tools to solve practical problems?

**Grenander:** Yes. In the past, some colleagues have accused me of writing in a purely mathematical style. That may be true, but the original problem under investigation often arose from real life. When confronted with technical difficulties, I often took a step back to look at a problem under more generality, not for the sake of generality, but to simplify the problem itself.

## VISITING BROWN UNIVERSITY (1957–1958): REKINDLING INTEREST IN COMPUTERS

**Mukhopadhyay:** Subsequently, you visited Brown University, did you not?

**Grenander:** I received a letter from Brown University asking whether I would like to come for a year. The invitation came from William Prager. He liked my work with Szegö and perhaps saw some connection between this piece of work with some of his own research. Prager was a very interesting German professor of applied mathematics. He created the field called plasticity theory within continuum mechanics and he had formed a wonderful group at Brown. This type of research program was initiated by the Navy at the end of World War II at Brown and at one other place that later became known as the Courant Institute. I decided to visit Brown as a professor of probability and statistics for nearly a year and a half in 1957–1958.

**Mukhopadhyay:** How was life at Brown in 1957?

**Grenander:** Of course, Brown was a small school, but it had a number of very strong departments. The applied mathematics group was one of the strongest. The applied mathematics program is still one of the very best. In 1957, the program was totally dominated by research in both fluid and solid mechanics, so I was an outsider from day one! I think that Prager's idea was to create a broader horizon.

**Mukhopadhyay:** Did you have any teaching obligation at Brown?

**Grenander:** I was teaching both probability and statistics. I was also collaborating with Walter Freiberger, who was the head of the computing lab: he was working with an IBM computer. Mr. Watson, a Brown alumnus, continued donating IBM computers to this lab. I was immediately drawn to computers and computing. Actually, I had been interested in computers long before I arrived here.

**Mukhopadhyay:** At what point did you feel initial attractions to work with computers?

**Grenander:** When I got out of the Army, I was married to Paj and I had to financially support us. To make ends meet, I became a teaching assistant under a scientist, Conny Palm, at the Royal Institute of Technology. Palm was a telephone engineer and he was a very interesting man. Early on, he gathered 10,000 used telephone relays, coupled them together, and built a computer. It was a relay computer, not an electronic computer. Such machines had been built at Harvard earlier. This was an interesting machine and it did not cost anything. The first time I used it for computing in time series, it had to be programmed with cables! Palm wrote a Ph.D. thesis in artificial teletrafficking (that now goes by the name queueing theory), which was published in German. I still have that monograph with me.

I carried that book in my rucksack when I joined the Army. It was a beautiful thesis. Well, I liked it, but no really influential people I knew at the time liked this work much. Later, D. G. Kendall was the one who publicized the real importance and usefulness of this highly original machinery called *Palm theory*. Now of course, Palm theory is very well regarded.

## BACK TO STOCKHOLM (1958–1966)

**Mukhopadhyay:** After your visit to Brown, did you return to Stockholm?

**Grenander:** Yes, I went back to Stockholm in 1958. In 1959, Cramér became the Chancellor and he relinquished his chair in the institute. I became his successor at this institute, where incidentally insurance mathematics came before mathematical statistics.

**Mukhopadhyay:** During those Stockholm days, did you supervise any thesis work?

**Grenander:** Yes, I had students working on various projects. Most of them wrote Licentiate theses and only a handful wrote Ph.D. theses. There was one very interesting thesis on insurance mathematics by Harald Bohman. There was another extremely good thesis written by Matérn. He was not my student really. I mean, I signed his documents because I was the professor, but actually he did all the work completely independently on his own. He was working in forestry. In a forest, trees grow densely in one part and not so densely in another part. Matérn had the brand new idea to describe the variation with the help of stochastic processes. Then he estimated, for example, the total volume of timber in a forest. This was a very important practical problem and Matérn's work was very fundamental in forestry research. This work had subsequently opened up the field of spatial statistics, made enormously popular and brought to the center of attention by Mardia, Besag and others.

Matérn was very original and clever. He learned all the mathematics on his own and I think that the Bayesian idea was his own too. He wrote his Licentiate thesis in Swedish as I did and hence not very many people read this work at the time. His Ph.D. thesis was in English though. Matérn, however, was severely criticized for his approach, because at that time one was not supposed to do Bayesian work! (Laughs)

**Mukhopadhyay:** Any other notable students?

**Grenander:** Sven Erlander's work in statistical problems led to optimization in general. He became a leader in the field, later becoming the President of Linköping University. He has recently retired. Per Martin-Löf wrote a Licentiate thesis on probability measures on semigroups. It was an extremely nice piece of work. There were others, for example, Walldin and Ekman, who had worked on operations research in banking.

## VISIT TO IBM WATSON LAB (1966): ENCOUNTERING PATTERN RECOGNITION

**Mukhopadhyay:** Please tell me about your visit to IBM.

**Grenander:** In 1966, I got an opportunity to spend one-half year at IBM's Watson Lab in Yorktown Heights, New York. This was a wonderful experience as I tried to learn computer programming and other related things.

The chairman of the mathematics department once told me to look into a field called pattern recognition. I heard about this field earlier, but did not pay much attention to it. I was given a list of nearly 500 publications to walk through!

**Mukhopadhyay:** What was your impression?

**Grenander:** I read many articles, perhaps 100 of them. I found them extremely boring! The treatment was mathematically trivial and the approach was often impractical. The literature at the time did not really handle "patterns." The field was growing fast and there were some successes, but these were very few and far between.

**Mukhopadhyay:** Who were some of the major contributors in pattern recognition as you were making your way?

**Grenander:** K. S. Fu was one of them and he was certainly one of the better ones. He wrote several books and important papers.

## PERMANENT MOVE TO BROWN UNIVERSITY (1966)

**Mukhopadhyay:** When did you join Brown University permanently?

**Grenander:** I came to Brown permanently in 1966 right after visiting IBM and I have stayed here for a very long time.

**Mukhopadhyay:** Did Harald Cramér visit Brown?

**Grenander:** Yes, he did. You know that Cramér came to visit the University of Connecticut in the early 1980s.

**Mukhopadhyay:** Oh yes. He visited the University of Connecticut in 1980 for the videotaping under the American Statistical Association's (ASA's) Filming of Distinguished Statisticians series. Cramér's lecture is preserved in ASA's archive.

**Grenander:** Cramér came to Brown on his way from the University of Connecticut and gave a series of three lectures. I remember the enthusiasm shown by our graduate students at the time when they heard him deliver rather old-fashioned and yet very elegant lectures. At that age, he was hard of hearing, but he was always a giving, caring and accommodating person. I loved him dearly.

When Cramér turned 85 or so, he started one of his last projects to learn ancient Greek and he managed to read Homer's *Odyssey* in Greek. You know that he made a trip to India when he was quite old.

**Mukhopadhyay:** Yes, I saw Cramér at the Indian Statistical Institute–Calcutta in December, 1974. Ulf, you talked about Feller before. Would you mention some of Feller's activities at Brown?

**Grenander:** Brown University housed a department of "history of mathematics" that was a small but very good department. Its leader was Neugebauer, a famous person who came from Germany. Feller had been the editor of *Mathematische Zentralblatt*, a sort of *Mathematical Reviews* from Germany. He arrived at Brown with a complete list of active mathematicians, including all the important names from all over the world. The number of active mathematicians at the time may be of some interest, and it was 300! Aided by this list, Neugebauer initiated editing the *Mathematical Reviews* from Brown. Later, Feller took over the editorial responsibilities and continued in that capacity for years until he left Brown and moved to Cornell.

**Mukhopadhyay:** What was Feller's personality like?

**Grenander:** He was a very colorful person. I remember that he once explained at a party how World War I really started. His view was quite different from the conventional one. Feller was very entertaining. He was a great storyteller and he told stories about practically anything imaginable.

**Ph.D. Student Advising**

**Mukhopadhyay:** Would you mention some of your graduate advisees from Brown?

**Grenander:** I had Rick Vitale and Don McClure as my students. Rick, of course, is your colleague at the University of Connecticut. Don McClure is with the applied mathematics group at Brown. Don took up the position to head the department for two terms and became my boss. He treated me well!

**Mukhopadhyay:** Ulf, were you pushing the frontiers of pattern theory?

**Grenander:** Yes, much work in pattern theory was done in the 1970s and 1980s by my students and me. Prior to my arrival at Brown, I wrote a little book about probabilities on algebraic structures to study probability distributions on groups, semigroups and things like those.

Chii Ruey Hwang was one of the best and wrote a thesis on limit theorems in pattern theory. Also, Yun-Shyong Chow was a brilliant mathematician. Both Chii Ruey and Yun-Shyong were very theoretical, but I wanted them to appreciate and work in applied settings too. That was, however, not to be, because neither enjoyed working with computers! Jointly with Yun-Shyong and Dan Keenan, I wrote the book *HANDS*: *A Pattern-Theoretic Study of Biological Shapes* (Grenander, Chow and Keenan, 1991).

**Creation of the Swedish Institute of Applied Mathematics**

**Mukhopadhyay:** At some point, did you not shuttle between Providence and Stockholm a few times a year? What was this for?

**Grenander:** Let me go back to something that I did when I left for Stockholm in 1958. I encouraged my graduate students to go out in the real world and find problems whose solutions required mathematics. For example, I had two students, Ekman and Walldin, who worked jointly and I asked them to look into problems in banking.

**Mukhopadhyay:** How did you proceed?

**Grenander:** I wanted to initiate a Swedish institute in applied mathematics, where graduate students would be paid to spend a couple of years working on practical problems from the industry and government. Then students would write theses about the problems and defend them. We needed a good deal of money for this operation. I tried to convince the governmental authorities in Sweden, but I received absolutely no response. Eventually, I moved to Brown in 1966. After I left Stockholm, I sensed some movement and interest on the part of the Swedish establishment. At that time, they asked me to return to Stockholm and initiate the institute. I was not sure if I should do so permanently, as I had already fallen in love with Brown.

**Mukhopadhyay:** But that is not to say that you were not interested to start that institute.

**Grenander:** Right. Deep down, I was always seriously interested in building the Swedish Institute of Applied Mathematics. Initially, I was commuting between Providence and Stockholm a number of times per year. I recall that this arrangement continued for four or five years and it was exhausting! The Institute finally took shape. Daniel Sundstrom, Goran Borg and Germind Dahlqvist helped to get it off the ground. I discovered later that the typical indifference and resistance from the government against starting such institutes were actually less under the Swedish Prime Minister, Mr. Palme.

**Mukhopadhyay:** This was a major interdisciplinary initiative. Who took care of its day-to-day operation?

**Grenander:** For a while (1971–1973), I was the Scientific Director of the Institute even though physically I was at Brown, which was little awkward, but we had a group of scientists on the board including Dahlqvist. Students and scientists from the Institute had wonderful opportunities to work on hundreds of practical problems.

### Teaching and Mentoring Undergraduate Students

**Mukhopadhyay:** Ulf, for a number of years, you have been teaching only undergraduate students at Brown. What is your motivation?

**Grenander:** The undergraduate students at Brown are very hardworking. I would say that our undergraduates are motivated to work harder than their peers. They pay so much! (Laughs)

From the very beginning, I found many undergraduates at Brown really interesting. I liked them so much that after being here for 20 years or so, I decided to teach only undergraduate courses.

**Mukhopadhyay:** Are you talking about required undergraduate courses or seminars?

**Grenander:** I have taught small senior seminars with five or ten students. I aim at exposing the students to ideas of mathematical experiments. I would pose a problem and first ask the students to guess the solution rather than solve the problem outright. At that point, I would ask students to use computers to evaluate all sorts of things and then examine the output themselves (see Grenander, 1982). I encourage students to form their own hypotheses and somehow approach to prove or disprove a hypothesis. These exercises often turn out very well and this involvement is satisfying too.

**Mukhopadhyay:** Any specific example?

**Grenander:** In my exclusive encounters with a select group of undergraduate students over the past ten years, I have met a number of very bright students. Once we worked on a big project and we had so much fun. It went over two semesters I recall. We had 20 students, divided into three groups, simulating three competing insurance companies. I wanted them to develop computer programs that would simulate successful behavior of insurance executives. Undergraduate students started working on the project, but eventually this led to a very fine Ph.D. thesis for Pat Burke, who is now associated with a university in Ireland.

### The Pattern Theory Group

**Mukhopadhyay:** When and how did you initiate the pattern theory group at Brown?

**Grenander:** A couple of years after 1970, I started organizing a pattern theory group at Brown. You will recall that I used to work with Walter Freiberger, who by then returned to statistics. Then Don McClure joined the group. A few years later, we had a young visitor, Stuart Geman. Stuart first came to Brown for a brief visit, I think. Do you know him?

**Mukhopadhyay:** Yes, I have met Stu Geman a few times.

**Grenander:** When Stuart visited Brown, I vaguely recall that he had a job somewhere else. He was just out of MIT with a Ph.D. thesis written under Herman Chernoff and instantly became interested in what we were doing at the time. He actually left the other job opportunity and decided to join our group. That was a great recruit for the pattern theory group.

**Mukhopadhyay:** Stu Geman surely attracted other people to join your group.

**Grenander:** Yes, he has been a tremendous help. Stuart attracted others including Elie Bienenstock from France and, more recently, David Mumford from Harvard. Both Elie and David have been with the group for a number of years. I should also mention Basilis Gidas, who came to us around the mid 1980s. He had already made a career in the field of partial differential equations; then he switched to physics and later he decided to join the pattern theory group. Now there is a program at Brown that goes by the name brain science. I think that it is one of the best and most active programs in the country, including pattern theory models and knowledge representation. It is not the biggest program, but certainly it is one of the best.

**Mukhopadhyay:** Within this group, initially how did you exchange ideas?

**Grenander:** Around 1980, we started an informal seminar series to discuss problems and projects. We did not necessarily have formal presentations. There were not that many of us at the time. This forum was dedicated to exchange of ideas and to inform what others in the group were doing. That used to be lot of fun.

## PATTERN THEORY AND BRAIN SCIENCE

**Mukhopadhyay:** You are a pioneer in pattern theory. To a layman, what is pattern theory?

**Grenander:** Let me start by saying that the objective is the same as in any other areas of science, namely we try to understand our surroundings, but the emphasis in pattern theory is on the actual *act of* knowledge and *act of* understanding. The key phrase is "act of," that is, we want to learn the process of understanding and the emphasis is not necessarily on specifics of what it is that we understand. Pattern theory is more like mathematics of knowledge representation (see Grenander, 1996).

**Mukhopadhyay:** Is this area related to some of the things that you had developed before?

**Grenander:** I realized that any pattern should be described with the help of algebraic objects and I was drawn to my earlier research on probabilities on algebraic structures. My goal was to put appropriate probability measures on the resulting regular structures.

I refer to pattern theory as the intellectual adventure of my life. I have applied this theory to solve a number of practical problems in biology and medicine. I am happy to say that I have done some of these works jointly with my buddy, Michael I. Miller, a prolific scientist. He is rich with ideas and he is lots of fun to work with. It has been great to collaborate with Michael in preparing the book, *Pattern Theory*: *From Representation to Inference* (Grenander and Miller, 2005).

In 1993, I published my chef d'oevre, *General Pattern Theory*, with Oxford University Press. It includes some definitive results as well as many more tentative ones. I hope that the younger generations of mathematicians will take up some of these topics.

**Mukhopadhyay:** Did you interact with David Mumford?

**Grenander:** Yes. I was greatly influenced by David Mumford. In a lecture he once said, "You should go out and measure the world." In pattern recognition, researchers have been trying to accomplish just that for over 50 years, but not so successfully though. New ideas and new directions are needed.

**Mukhopadhyay:** Have you approached such problems in a different way?

**Grenander:** Yes, I thought that there should be an analytical model first. One may refer to Grenander and Miller (1994) for specific details.

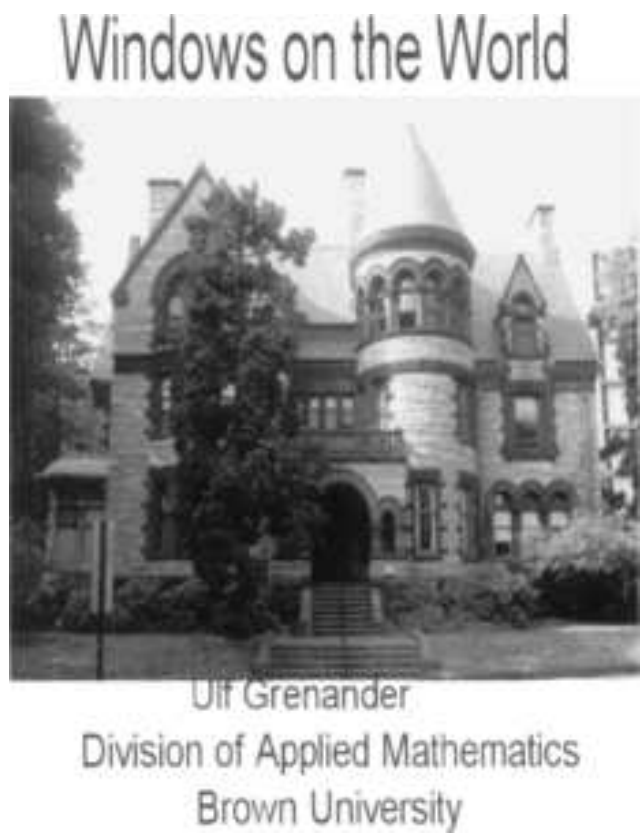

FIG. 8. *The cover for one of Ulf Grenander's CD-ROMs.*

Mumford and I studied the mathematics of natural scenes such as a forest or cityscapes. Everywhere we saw the characteristic cusp. I did not realize at first that this was not simply a qualitative agreement. It turned out to be much more important than that. My colleague, Anuj Srivastava from Florida State, and I collaborated on several things and bumped into a universal law! This is an extremely promising approach in pattern theory now.

**Mukhopadhyay:** You also have interest in what you call mind theory. What is it briefly?

**Grenander:** Several people in our pattern theory group have worked on the brain's physiology and neural system. They actually measure things microscopically inside the brain. In the future, one will gain better understanding of how the brain's neural system functions on a global level. Then one may be able to mathematically describe human thought. In this discussion, I do not include simply the logical thoughts alone. One should be describing, for example, human emotions, doubts and an act of lying or talking. That will surely happen some day, but my present approach is completely speculative, based on introspection rather than hard data.

## STATISTICAL GROOMING

**Mukhopadhyay:** Were you influenced by some of the mathematicians and statisticians from England and America in shaping your career?

**Grenander:** Sweden was surrounded by the German Army during WWII, so we lived in complete isolation and did not get any mail from other countries. During my upbringing, I did not see any book or journal or newspaper from other countries except Germany. In high school, I used to buy only German books in mathematics and some of them were very good. Even in college, German influence was dominating. Years later, I found scientific materials arriving in Sweden from America, England and other places.

**Mukhopadhyay:** You surely learned some basics from the works of K. Pearson, R. A. Fisher and J. Neyman, did you not?

**Grenander:** Of course, I did very much. My fellow students and I read a series of papers written by J. Neyman and E. S. Pearson. Those papers were very intense and filled with brand new ideas. Their approaches to statistical inference were revolutionary at the time. Neyman and Pearson cleared up lots of the mess within statistical theory that existed before they came along.

I also learned from the writings of some people whose names may not sound very familiar today. I learned from Tschuprov's correlation theory. The original work was in Russian, but I read it in German. I read Elderton and we had to memorize his eight classes of distributions! I also read Whittaker and Robinson's work. I heard about fast forward Fourier transforms before they were invented! There was a book that schematically gave the way to calculate Fourier sums!

**Mukhopadhyay:** The formal logic was never foreign to you, but when you first encountered the works of Fisher or Neyman, did the deductive nature of logic in statistics bother you?

**Grenander:** From the point of view of logic, I did not feel disturbed. Statistics is not formal logic, but it develops machinery to handle uncertainty. I loved that aspect, I really did. The world is made up of uncertainties. There are only a few certain things! The field opened up wonderful opportunities for me. When I learned statistical mechanics for the first time, I thought that it was absolutely wonderful.

I was delighted to enter the field of statistics and probability, and yet I have stayed close to mathematics. I sensed a great future in this and I am very happy that I became involved.

**Mukhopadhyay:** Did your visit to Chicago influence your statistical views in any way?

**Grenander:** In Chicago, I experienced the wonderful cultural milieu. That was incredibly invigorating for my mind and soul. My statistical ideas, however, were not too influenced by that environment at the time. It was fashionable to formulate every statistical problem as one involving decision functions. We still do some of that, but at the time this statistical culture was taken to its extreme! I found that approach repellent, too abstract, too general and with too little substance. Often a very general approach does not lead to many fruitful clues!

**Mukhopadhyay:** A complete class, for example, may be too large to handle!

**Grenander:** Nitis, you are absolutely correct. Instead of going after the most general theorem, staying mindful about specifics may sometimes prove helpful. I am a believer of that. The earliest influence on my mathematical views came from Beurling. I mentioned that before. He was the most powerful intellectual I ever met.

**Mukhopadhyay:** Some of your own writings have been much too abstract, I might add.

**Grenander:** Yes, some of the things I have written over the years are a little abstract I suppose. I realize that. I once wrote a book called *Abstract Inference*. This was probably not a catchy title for a book! Later someone pointed out that what I had discussed there should not have been called *abstract* because the approach I proposed was rather *concrete*!

## RADICAL VIEWS ABOUT STATISTICS: PARAMETRIC VERSUS NONPARAMETRIC

**Mukhopadhyay:** Ulf, would you say that your view of statistics has changed significantly?

**Grenander:** Yes, I would say that. I have certainly changed my view about what statistics ought to be, but this has evolved over a long period of time.

Early on, I used to do lot of research in medical statistics involving bioassays, clinical trials and pharmaceutical problems. So I am also guilty of having committed hundreds of analyses of variance! Later on, I have wondered whether these were the right things to do. Under the assumption of normality in each situation, we obtained what was called an *exact test*! There was nothing *exact* about those tests, but these were still called exact tests. The name was a misnomer because nobody believed in the assumptions that were made in the first place. I have come to the conclusion that those methods should be used with plenty of caution.

**Mukhopadhyay:** You are not asking others to abandon small sample exact tests, are you?

**Grenander:** Small sample exact tests should be taught in classrooms if for no other reason than that these tests are being used extensively in applications. This material should be appropriate to teach in a course if the treatment is accompanied with a good dose of skepticism.

**Mukhopadhyay:** What will you propose to a practitioner as an alternative approach?

**Grenander:** I think that a nonparametric approach is often the way to go, because this does not involve as many assumptions. I also really believe in large sample theory. I have come to believe in these kinds of methodologies.

**Mukhopadhyay:** But there may be occasions where one would have fairly valid reasons to work with a parametric model.

**Grenander:** I must admit that there may be occasions where one would have valid reasons to postulate particular parametric models. Let me give you my favorite examples. In signal processing, for the kinds of things Stuart Rice did so successfully, the customary assumptions could be validated and parametric formulations worked. The same may become true about two- or three-dimensional image processing. We know that sometimes there are valid underlying models. In statistical mechanics, there are valid parametric models and they work. This is one of the great success stories in science! The field of insurance mathematics is similar in this sense.

But I am not inclined to interpret the significance levels as probabilities. They may be useful for calibration purposes or for setting up some standards for comparison. Persi Diaconis has expressed similar views better than I have.

**Mukhopadhyay:** Is it fair to say that because in much of your work you adopted Bayesian approaches, you cannot accept significance levels as probabilities?

**Grenander:** No, the reason is not necessarily a reflection of my inclination to often use Bayesian approaches at all. My feeling was rather shaped by the fact that the assumptions behind those *exact tests* were so stringent.

## DOOB, WIENER, HÁJEK AND DALENIUS

**Mukhopadhyay:** Given your lifelong contributions in time series and stochastic processes, will you please remark briefly about Joe Doob?

**Grenander:** I met Joe Doob a number of times. I remember in particular that once I gave a seminar at Urbana–Champaign, perhaps in 1958. I was in Chicago at the time and some of us drove down to Urbana. I gave a lecture about optimal regression in a stationary stochastic process or something like that. I heard later that Doob gave his own version in a seminar the following week and the topic was "What Grenander really meant." (Laughs)

**Mukhopadhyay:** Were you influenced by Doob's work?

**Grenander:** Of course Doob's work influenced my research in time series and stochastic processes very much. If you ask my wife, Paj, she will testify to how much I really suffered through one of Doob's first papers on stochastic processes. He had a lemma that I thought was completely obvious, but, I also thought to myself, "I could not possibly be right, because if it was indeed as simple as I thought it was, then why did Doob supply its proof!" It took me weeks to understand why that one particular lemma was not really obvious! (Laughs)

**Mukhopadhyay:** Have your research interests in stochastic processes changed?

**Grenander:** At the very beginning of my career, I was very much influenced by Doob's and Kolmogorov's work on stochastic processes. Then, of course, my own teachers influenced my thoughts. These two major and yet different kinds of influences empowered me to create the field of statistical inference in stochastic processes. That was my first sincere interest in science.

**Mukhopadhyay:** Do you recall interactions with Norbert Wiener?

**Grenander:** Norbert Wiener visited Stockholm repeatedly when I was there. We used to go to lunch together and I remember that he asked me the same question again and again, "Why am I not appreciated as much in my country as Kolmogorov is respected in the Soviet Union?" We must have gone through this conversation at least three times! He genuinely felt that he was not appreciated in America and it was clear that this made him feel unwanted, but, of course, he was appreciated a lot in America, perhaps more so than anybody else! He was treated like a God but he was not satisfied with the level of recognition that he had received.

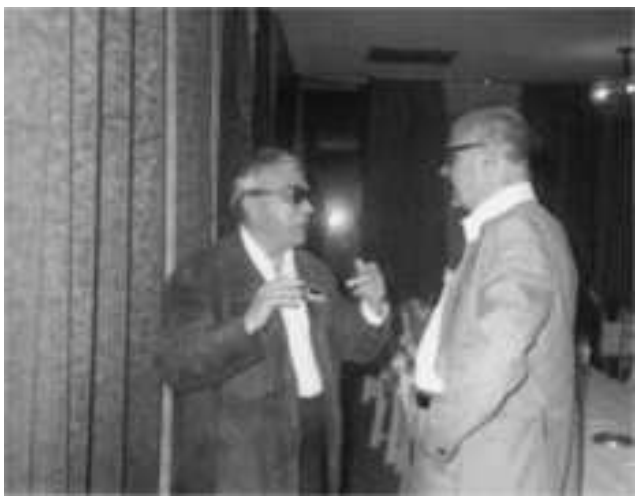

Fig. 9. *S. Kh. Sirazhdinov (left) and Ulf Grenander.*

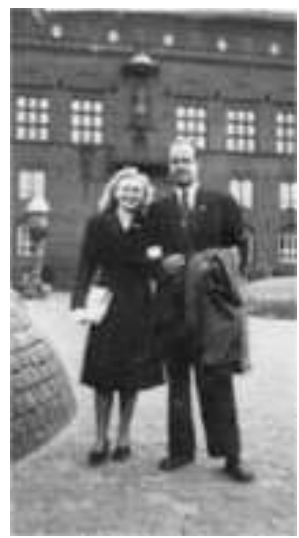

Fig. 10. *Ulf Grenander with fiancée, Paj, at the Copenhagen Math Congress, 1946.*

**Mukhopadhyay:** At any time did your journey cross the path of Hájek?

**Grenander:** Hájek came to Stockholm and stayed there for quite a while. He did interact with both Cramér and me. Tore Dalenius and Hájek had common interests (Hajek, 1964) and so they interacted with each other too. Tore spent almost 20 years at Brown as a visiting professor. Unfortunately, he passed away recently.

**Mukhopadhyay:** I am extremely sorry to hear that.

**Grenander:** Perhaps Tore and Hájek were working on something related to optimal stratification, one of Tore's main contributions in statistics. Those were some very beautiful pieces of work (Dalenius, 1950; Dalenius and Hodges, 1957, 1959).

Tore and I were very close friends. We spent lot of time together discussing many things. Although we never wrote a paper together, we had many scientific interactions and exchanges. Tore did not have any formal mathematical training, but the mathematics he knew, he used to a great advantage. Should I tell you a story?

**Mukhopadhyay:** Ulf, of course. Please go right ahead.

**Grenander:** My wife, Paj, used to subscribe to the *Swedish Medical Journal*. I think that I read it more than she did. In the 1950s or 1960s, I saw an awful paper there advocating computer usage in medicine. I am all in favor of computers and I love them, but the treatment in that paper was so naive! I showed this work to Tore and said, "Look at this incompetence. Perhaps this piece is published as a joke." Tore instantly replied, "No, no, you can never underestimate incompetence."

**Mukhopadhyay:** (Laughs) I did not realize that Tore had such sense of humor.

**Grenander:** But, Nitis, listen, that was not the end of the story. Tore said, "Let us try something different." He wrote a paper as a joke and it pretended to support the finding of the other published paper and said that it was a "wonderful paper"! He sent off this paper to the same journal for publication and it was accepted! Tore's paper came out in print in the *Swedish Medical Journal*! The authors of that earlier paper became furious after realizing that Tore made fun at their expense. We had rather unpleasant discussions with them later.

## MEMORABLE CONFERENCES

**Mukhopadhyay:** Is there a conference that is especially etched in your mind?

**Grenander:** I particularly remember one conference. I was still in the Army. I managed to get a leave of absence and appeared at a conference in Copenhagen at the end of the War in my uniform with my girlfriend, Emma-Stina, delightfully called Paj. This was my first scientific Congress and it was the best conference I have ever been to. I heard many wonderful lectures there. There was one special talk that I had not anticipated and it was given by an American mathematician, Henry Wallman.

**Mukhopadhyay:** What was it about Wallman's lecture that impressed you?

**Grenander:** I thought that he was one of the most abstract mathematicians possible. He was writing a book on something called dimension theory. Interestingly, he was also a member of the American Communist Party. At the end of the War, he got into trouble in America. I do not know if he had to leave America or not, but he ended up in Sweden. While living in Sweden, he suddenly switched his profession and became an electronic engineer. Wallman had an extremely mathematical as well as a practical bent of mind. At the Copenhagen Congress, he gave a talk about computers. This was really precomputer days! Wallman mentioned building an internal computer library to enable evaluations of some of the usual functions such as sin, cos, exp and so on. What a sci-fi idea at the time! I was immediately turned on.

After the Congress, Paj and I went to visit Elsinor Castle in Northern Denmark, where Hamlet's characters had supposedly lived. There I kept on talking about Wallman's wonderful idea for hours! Paj was surely bored to tear, but my enthusiasm was simply unstoppable! (Laughs).

**Mukhopadhyay:** How about conference trips to Greece?

**Grenander:** I first visited Greece in 1965 for a conference in Loutraki, some distance from Athens. That was the first time I talked about pattern theory in a conference! There was no such field at the time and I suggested that there ought to be one.

**Mukhopadhyay:** Did you not visit France and Italy?

**Grenander:** Oh yes, I visited France many times. I started visiting France very early in my career. Robert Fortet from the University of Paris and I had many special common interests including projects on inference in stochastic processes. I had contacts with one of his students, Édith Mourier, from France. Mourier wrote a very interesting thesis about probability measures on a Banach space. There were other very good probabilists, for example, Meyer.

I had many more visits to Italy though. At one point, I was visiting Rome a couple of times a year. In the late 1980s and early 1990s, I worked with a group from Rome that was led by Mauro Piccioni, who had done work on inference in stochastic processes. That brought me closer to this Italian group. However, ultimately we did not work together on inference in stochastic processes though. Instead, we collaborated in the field of image processing. I worked on a pattern theoretic setup having a structure of image algebra together with Yali Amit, who spent some time at Brown.

## COLLEAGUES AND FRIENDS

**Mukhopadhyay:** Did you meet P. V. Sukhatme? Any recollections?

**Grenander:** P. V. Sukhatme visited Stockholm once that I am aware of and he talked about sampling. He had his own ideas and approach to gather large agricultural data, but his methodology came under fire in the light of a more acceptable approach due to P. C. Mahalanobis. Neither compromised individual scientific ideals.

**Mukhopadhyay:** Did you not collaborate with J. Sethuraman?

**Grenander:** Yes, I did. At some point, I spent some time in Florida. I met Sethuraman there and we worked on limit theorems in pattern theory. He did some beautiful stuff and some of this was published in the *Journal of Multivariate Analysis* and the *Journal of Theoretical Probability* (Kurien and Sethuraman, 1993a, b). Later, Sethuraman and I worked on something called *conformation of molecules* (Grenander and Sethuraman, 1994). We adopted a pattern theoretic approach to describe how two geometrical figures merge in mathematical terms. By that time, I started moving away from theory to more practical applications of pattern theory.

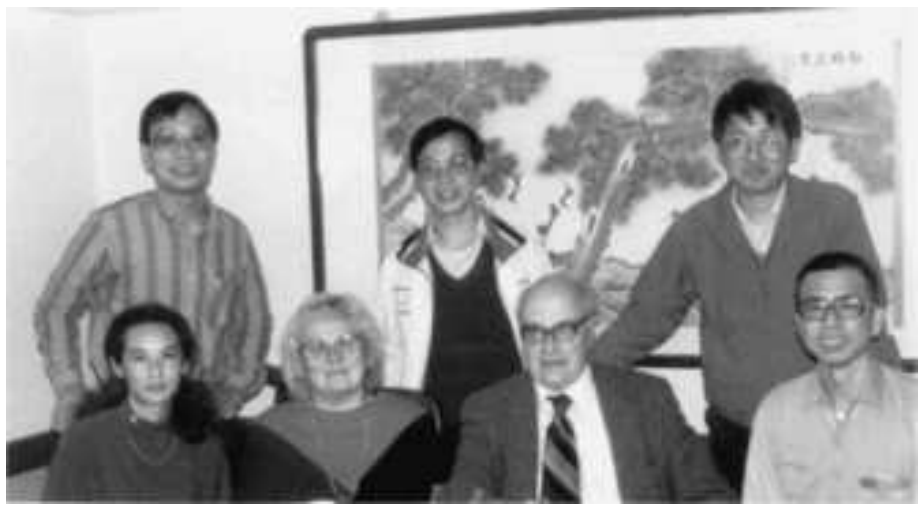

Fig. 11. *Ulf Grenander with his wife, Paj, and Chinese colleagues.*

**Mukhopadhyay:** Do you have anything on C. R. Rao?

**Grenander:** I have met C. R. Rao many times. During a recent visit to Penn State, I gave a talk where I pointed out that I used a mathematical technique borrowed from him. The first time I heard about it was during one of Cramér's seminars many years earlier. I suppose that I must have been a student at the time. Cramér discussed a result of C. R. Rao about optimal estimation. Cramér said, "It is a beautiful result, but can it be true? It is too beautiful!" Of course, the result turned out to be true and it has come to be known as the Rao–Blackwell theorem, but during 1946–1947 when I was a student, it was a very surprising result. Apparently, it even surprised Cramér.

**Mukhopadhyay:** Do you wish to mention M. M. Rao?

**Grenander:** I was very happy when he contacted me a few years ago and mentioned that he was writing a book partly based on my Ph.D. thesis (Grenander, 1950), but it would be much more advanced and fully developed.

**Mukhopadhyay:** Your friendship with Kanti Mardia goes way back. Any thoughts?

**Grenander:** Kanti Mardia and I have been very good friends for many years. I have visited Leeds and we have met in this country on different occasions. My research interests are very similar to Mardia's as you know.

Kanti has really built a very fine group at Leeds and I believe that the group is broadening the scope of research on many fronts. The Leeds Applied Statistics Research Workshops have turned out to be quite some international events. The credit goes to Kanti and his group.

**Mukhopadhyay:** Any recollections about Geoff Watson?

**Grenander:** I had some interactions with Geoff Watson about statistics on spherical data. I had a graduate student, Bjorn Ajne, at Stockholm. Bjorn and I used to meet Geoff in all sorts of different places to discuss research problems and exchange ideas. At that time, Geoff was at Johns Hopkins in Baltimore, but later he moved to Princeton.

## GENOME, DATA MINING AND BEYOND: FUTURE PREDICTIONS

**Mukhopadhyay:** Would you speculate about some of the areas in statistics and probability that may be in the forefront of science in the next 20 or 30 years?

**Grenander:** Let me turn the question around and say this. The field of statistics has evolved from applications. Probability theory was originated and energized by problems arising from gambling! Important stuff in both probability theory and statistics was developed and nurtured by actuarial mathematics. The life tables that started in 1750 or so had a big influence in our field. Statistical mechanics played a significant role in the history of the development of probability theory and statistics. What I think is going to happen is that both biology and medicine will become much more mathematical. Thirty years ago in Stockholm, I had many friends from the medical field and I tried to get their attention. I used to tell them to be skeptical about ANOVA (analysis of variance). I carried out more ambitious analysis with mathematics instead. My friends from the medical community were not impressed. They used to say, "This is all quite nice, but now is not the time for a change. Wait a little."

Now when I meet their younger generation, I find that some of them know mathematics well, perhaps because of their familiarity with computers. Now, I can see the revolution coming. Major impacts in our field in the next 20 or 30 years will come from biology and medicine. I believe that. Of course, probability theory and statistics will continue to grow on its own too, but real important advances will come from outside of statistical science. We ain't seen nothing yet!

**Mukhopadhyay:** What is your impression about data mining?

**Grenander:** Some neural scientists at Brown place 25 electrodes in a monkey's brain now. One could record responses from a 25- or 100-dimensional time series. Data mining there is a widely open field. If I were young, I would jump and grab this opportunity, but I must caution and add that to be able to do such things really well, it is not enough to be a mathematician. One must also learn important aspects of science along with computers.

**Mukhopadhyay:** But it is a reality that many "new ideas" may not withstand the test of time.

**Grenander:** That is true, but some not-so-good ideas may lead to important ideas much later. The name cybernetics is not heard as much today as in the 1950s, but, the basic principles of cybernetics laid down by Wiener many years ago are more relevant now than ever. In a way, that concept is still alive and well, and even today its influence is felt.

## THE IMMEDIATE FAMILY: INTIMATE THOUGHTS

**Mukhopadhyay:** Please tell me about sailing in the Baltic.

**Grenander:** Earlier in life, I loved sailing. The Baltic is excellent for sailing. It is never as cold as one might think because the Gulf Stream heats the water. There are nearly 5,000 small islands, known as the archi-pelago, outside Vastervik where we have a house. Navigating around these islands is not very easy. There are probably only a few hundred houses on all the islands combined! Most islands are vacant. A hundred or so years ago, only fishermen lived on those islands. Now the vacationers have all those islands to themselves.

**Mukhopadhyay:** Ulf, when did you and Paj first meet?

**Grenander:** The first time the two of us met was when Paj was only five and I was seven.

**Mukhopadhyay:** When did you two get married?

**Grenander:** I was 18 and Paj was 16 when we met again in high school. We became engaged on Paj's high school graduation day in 1945. We got married in 1946.

A little later, Paj started attending medical school, the Caroline Institute in Stockholm, which was very demanding. One day Paj came home from school and said that in ten days she will have her first test. She would need to memorize 1500 Latin names for the pieces of bones in a human body. I told Paj, "That would be impossible to do. I am sure that they would not ask just that sort of thing," but that was exactly what the test was about! Life for her became even harder when our children were on their way because she was still in school.

**Mukhopadhyay:** When were your children born? Where are they now?

**Grenander:** Three children were successively born in 1951, 1955 and 1957. We had a boy first and the next two were girls. Our son, Sven, is with the Jet Propulsion Laboratory at Cal Tech in space research and he is married to a girl he met there. It is interesting to note that she is of Irish origin, but was brought up partly in Spain in a Catholic school. So there we have a born foreigner in the family. (Laughs)

Then, our daughter Angela is an M.D. married to another M.D. from Iran. The younger daughter, Charlotte, is a physical therapist and she married a professor of medicine at UCLA and the Medical Director of Orange County. Charlotte is a personal trainer for actors in Hollywood.

**Mukhopadhyay:** How about your grandchildren?

**Grenander:** Angela is a pediatrician and she has given us four grandchildren, Alexander (17), Ariana (15), Nikolas (13) and little Tatiana (8). Our other daughter Charlotte has two children, Annika (11) and Anders (5). Our son has no children.

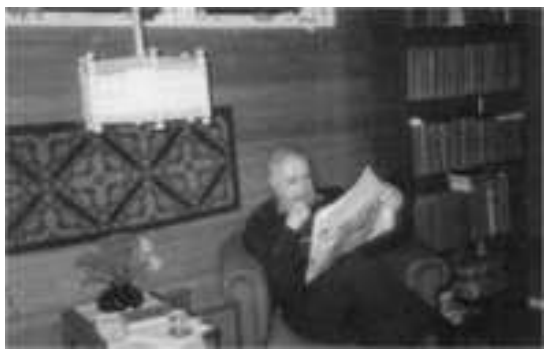

Fig. 12. *Paj Grenander's embroidery on the wall. Ulf Grenander in his summer house, Stockholm.*

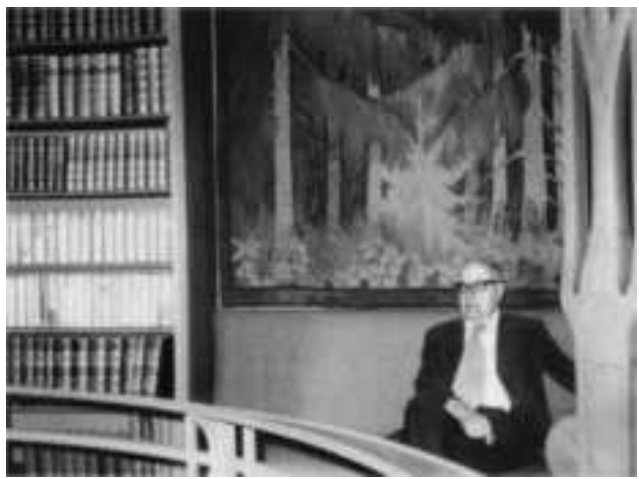

Fig. 14. *Ulf Grenander at Mittag-Leffler, 1982.*

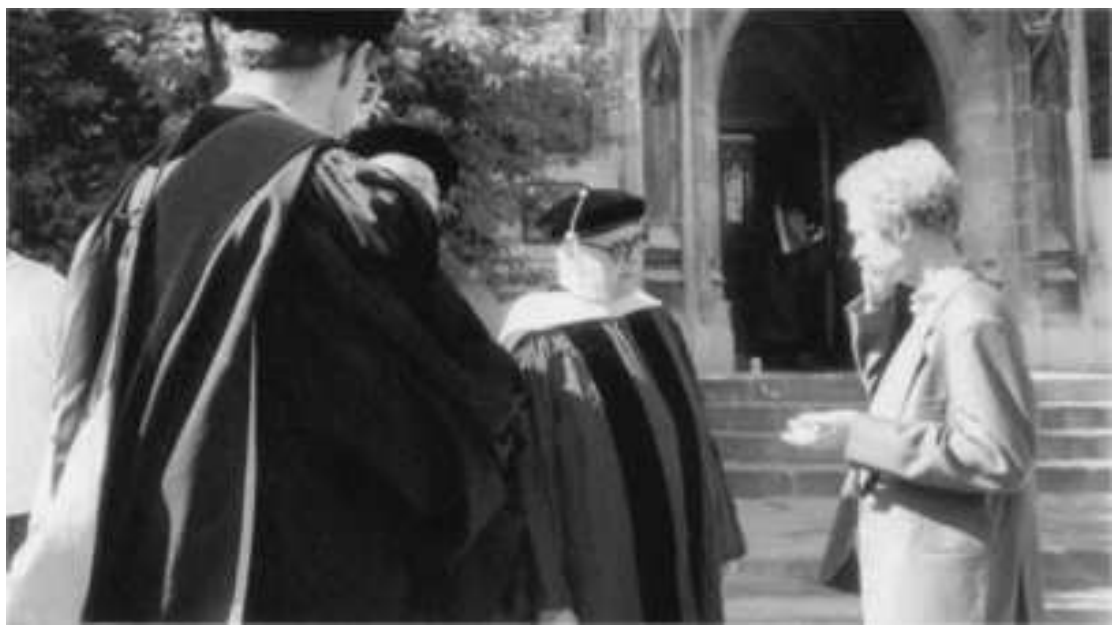

FIG. 13. *Honorary doctoral degree for Ulf Grenander (in the center) from the University of Chicago, 1993.*

I should mention that Alexander is my real buddy. He prides himself as very Swedish. His dad is of Iranian origin and Alexander told him one day, "Dad, you are a foreigner since you do not speak Swedish." (Laughs)

**Mukhopadhyay:** Every summer, your grandchildren go to Sweden to visit with you. That must be heavenly!

**Grenander:** Yes, indeed. Every year, Paj and I spend three summer months in Sweden, and a very important part of our interactions with the grandchildren unfolds when they come to visit Sweden. We live near water and boats would be waiting! So everyone has lots of fun.

Rufsan, our family pet, has been with us for a long time. She is very important to me and to my family. The book, *HANDS*, that I wrote with Chow and Keenan was not dedicated to Rufsan, but her favorite dog food was photographed to create some of the images used as illustrations.

**Mukhopadhyay:** Ulf, do you have a hobby?

**Grenander:** My main hobby is mathematics, but I also have some secondary hobbies including music, literature and history. History, especially modern history, interests me very much. I used to play bridge, but it took too much time away from mathematics, so I gave it up. But I must add that Paj has remained an enthusiastic bridge player.

**Mukhopadhyay:** I have glanced over Paj's colorful and wonderful embroidered pieces in your living room, family room and study. There is so much talent under one roof!

**Grenander:** Indeed, Paj is very artistic and she painstakingly creates fine artwork in typical Swedish style.

**Mukhopadhyay:** Ulf, I realize that you must feel proud of all the honors you have received. Would you please mention one or two that warm your heart and soul and perhaps make you think, "Yes, I made it."

**Grenander:** In 1965, I felt delighted when I became a member of the Royal Swedish Academy of Sciences. The academy hands out the Nobel Prizes in physics and chemistry, and it is decided by the members' votes. During ordinary meetings of the academy, each member receives a silver coin. I collected a number of them for my children and grandchildren. At the Nobel meeting, the King used to join and take part in the discussion, and each member received a gold coin.

**Mukhopadhyay:** You were awarded an Honorary D.Sc. degree from the University of Chicago in 1993. Please tell me about it.

**Grenander:** I received this honor in a special convocation, I suppose, because of my association with them in 1951–1952. I remember that I had to give a short speech after accepting the degree. The president of the university at the time was Dr. Sonnen-schein. So, I ended my speech by saying, "This University is in very good hands. It will be covered by a lot of sunshine." (Laughs)

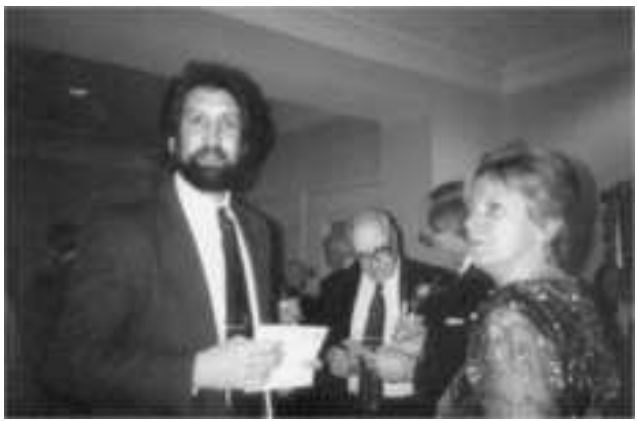

FIG. 15. *Michael Miller, Ulf Grenander and M.E. McClure.*

## LOVE FOR BROWN UNIVERSITY IN RETIREMENT

**Mukhopadhyay:** It is hard for me to think of Brown University and Ulf Grenander separately. In retirement, as you look back, what comes to your mind?

**Grenander:** I love mathematics and I have loved Brown University from day one. Brown not only gave me the splendid opportunity to challenge the frontiers of mathematics, it actually encouraged me. It pushed me to boldly pursue my research ideas, some new ways to do and apply mathematics! Brown University has been very appreciative of our work. I have to say that Paj and I have been really lucky and blessed. It has been a wonderful journey in life.

**Mukhopadhyay:** I heard that Brown appointed a new president recently.

**Grenander:** Yes, we have a new president, Dr. Simmons. The other day, she reiterated Brown's established philosophy that challenges both its faculty members and students to do daring things. These are the ideas and approaches perhaps not fashionable today, but they will break new ground in science in the future.

**Mukhopadhyay:** If you could go back and change anything, Ulf, what would that be?

**Grenander:** No, I do not think that I would change anything really. Perhaps I should have come to live in America earlier than when I did. I have gained so much from here in my intellectual life as well as my family life.

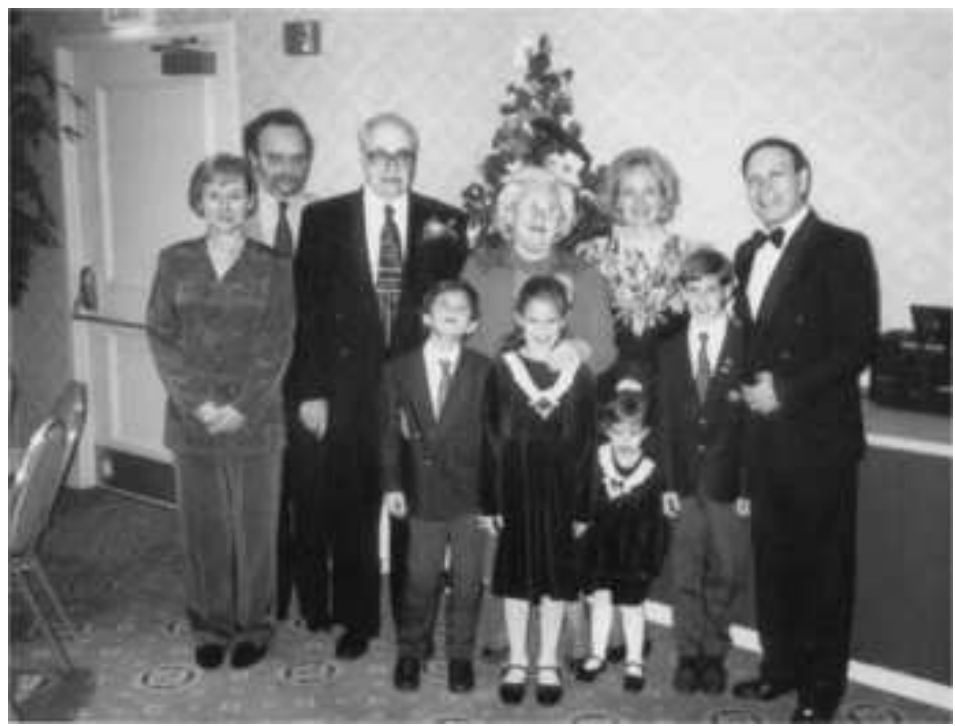

FIG. 16. *Ulf and Paj Grenander with family at their 50th wedding anniversary. Front row (from left to right): Nikolas, Ariana, Tatiana, and Alexander Raufi. Back row (from left to right): Nancy, Sven, Ulf and Paj Grenander, and Angela and Noori Raufi.*

**Mukhopadhyay:** By the way, for our readers, will you please mention your position now at Brown?

**Grenander:** Officially, I am a Research Professor here and what that means is that I am paid very little for working harder than I used to before retirement! (Laughs)

I have retired from teaching and in that sense I hold an emeritus position. I still spend a lot of my time working on pattern theory and mind theory.

**Mukhopadhyay:** Ulf, your energy is enviable. It has been a privilege to come to your home and have this conversation. I wish you and Paj a long, happy and productive life ahead with your wonderful family. Many thanks.

**Grenander:** Nitis, thanks to you as we close this conversation.

## ACKNOWLEDGMENT

I thank the Executive Editor and one of the editors for a number of helpful suggestions.

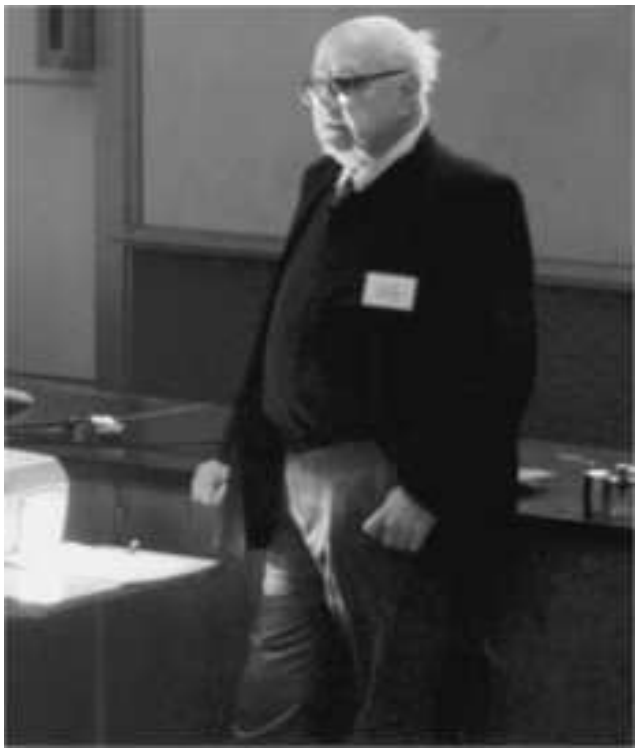

Fig. 17. *Delivering a special invited lecture at the University of Connecticut, Storrs, on April 6, 2002.*